\documentclass[10pt]{amsart}
\usepackage{amssymb,amstext,amsmath,amscd,amsthm,amsfonts,enumerate,latexsym,stmaryrd,multicol,geometry,graphicx}
\usepackage[usenames]{color}
\usepackage[all]{xy}
\newtheorem{thm}{Theorem}[section]
\newtheorem{lem}[thm]{Lemma}
\newtheorem{prop}[thm]{Proposition}
\newtheorem{cor}[thm]{Corollary}
\theoremstyle{definition}
\newtheorem{dfn}[thm]{Definition}

\newtheorem{setup}[thm]{Setup}
\newtheorem{rem}[thm]{Remark}

\newtheorem{conj}[thm]{Conjecture}
\newtheorem{ex}[thm]{Example}

\theoremstyle{remark}

\newtheorem*{claim*}{Claim}

\numberwithin{equation}{thm}
\def\cm{\operatorname{CM}}
\def\codepth{\operatorname{codepth}}
\def\D{\operatorname{D}}
\def\db{\operatorname{D^b}}
\def\depth{\operatorname{depth}}
\def\ds{\operatorname{D_{sg}}}
\def\e{\operatorname{e}}
\def\Ext{\operatorname{Ext}}
\def\ext{\operatorname{ext}}
\def\ge{\geqslant}
\def\grade{\operatorname{grade}}
\def\h{\operatorname{H}}
\def\Hom{\operatorname{Hom}}
\def\le{\leqslant}
\def\m{\mathfrak{m}}
\def\Max{\operatorname{Max}}
\def\Min{\operatorname{Min}}
\def\mod{\operatorname{mod}}
\def\nf{\operatorname{NF}}
\def\p{\mathfrak{p}}
\def\q{\mathfrak{q}}
\def\R{\mathbb{R}}
\def\reg{\operatorname{Reg}}
\def\rfd{\operatorname{Rfd}}
\def\sing{\operatorname{Sing}}
\def\spec{\operatorname{Spec}}
\def\syz{\Omega}
\def\T{\mathcal{T}}
\def\thick{\operatorname{thick}}
\def\V{\operatorname{V}}
\def\X{\mathcal{X}}
\def\xx{\boldsymbol{x}}
\def\Y{\mathcal{Y}}
\def\Z{\mathcal{Z}}
\def\ZZ{\mathbb{Z}}
\begin{document}
%\baselineskip 12pt
%\allowdisplaybreaks
\title[DHKK complexities for singularity categories and generation of syzygies]{Vanishing of DHKK complexities for singularity categories\\
and generation of syzygy modules}
\author{Tokuji Araya}
\address[Tokuji Araya]{Department of Applied Science, Faculty of Science, Okayama University of Science, Ridaicho, Kitaku, Okayama 700-0005, Japan}
\email{araya@ous.ac.jp}
\author{Kei-ichiro Iima}
\address[Kei-ichoro Iima]{Department of Liberal Studies, National Institute of Technology (KOSEN), Nara College, 22 Yata-cho, Yamatokoriyama, Nara 639-1080, Japan}
\email{iima@libe.nara-k.ac.jp}
\author{Ryo Takahashi}
\address[Ryo Takahashi]{Graduate School of Mathematics, Nagoya University, Furocho, Chikusaku, Nagoya 464-8602, Japan}
\email{takahashi@math.nagoya-u.ac.jp}
\urladdr{https://www.math.nagoya-u.ac.jp/~takahashi/}
\thanks{Ryo Takahashi was partly supported by JSPS Grant-in-Aid for Scientific Research 23K03070}
%\dedicatory{}
\begin{abstract}
Let $R$ be a commutative noetherian ring.
In this paper, we study, for the singularity category of $R$, the vanishing of the complexity $\delta_t(X,Y)$ in the sense of Dimitrov, Haiden, Katzarkov and Kontsevich.
We prove that the set of real numbers $t$ such that $\delta_t(X,Y)$ does not vanish is bounded in various cases.
We do it by building the high syzygy modules and maximal Cohen--Macaulay modules out of a single module only by taking direct summands and extensions.
\end{abstract}
\maketitle
%\tableofcontents
%%%%%%%%%%%%%%%%%%%%%%%%%%%%%%%%%%%%%%%%%%%
\section{Introduction}

Let $\T$ be a triangulated category.
Let $G$ be a {\em generator}, that is, an object of $\T$ whose thick closure equals to $\T$.
Let $X$ be any object of $\T$.
In their notable 2014 paper, Dimitrov, Haiden, Katzarkov and Kontsevich \cite{DHKK} define the {\em complexity}, which we call the {\em DHKK complexity}, of $X$ relative to $G$ as a certain real function
$$
\R\ni t\mapsto\delta_t(G,X)\in\R_{\ge0}.
$$
For the precise definition of a DHKK complexity, see Definition \ref{1}.
This notion is used in \cite{DHKK} to define the so-called {\em (categorical) entropy} of an exact endofunctor of $\T$.

Let $R$ be a commutative noetherian ring.
Denote by $\ds(R)$ the {\em singularity category} of $R$, which is defined as the Verdier quotient of the bounded derived category of finitely generated $R$-modules by perfect complexes.
Takahashi \cite{delta} presents the following conjecture.

\begin{conj}[Takahashi]
Let $G$ be a generator of $\ds(R)$.
Let $X$ be any object of $\ds(R)$.
One then has $\delta_t(G,X)=0$ for all nonzero real numbers $t$.
\end{conj}

The above conjecture holds when $R$ is a local hypersurface; this is an easy consequence of the fact that $\ds(R)$ is a $2$-periodic triangulated category; see Proposition \ref{25}.
Takahashi \cite{delta} proves the following theorem as the main result of the paper, which supports the conjecture.

\begin{thm}[Takahashi]\label{23}
Suppose that $(R,\m,k)$ is a local ring with an isolated singularity.
Let $G$ be a generator of $\ds(R)$.
Let $X$ be an object of $\ds(R)$.
Then the following assertions hold true.
\begin{enumerate}[\rm(1)]
\item
Put $c=\codepth R$ and $m=\max_{1\le i\le c}\{\dim_k\h_i(K^R)\}$.
Then $\delta_t(G,X)=0$ for all $|t|>\frac{1}{2}\log cm$.
\item
If $R$ is Gorenstein and $k$ is infinite, then $\delta_t(G,X)=0$ for all $|t|>\log(\e(R)-1)$.
\item
If $R$ is a complete intersection, then $\delta_t(G,X)=0$ for all nonzero real numbers $t$.
\end{enumerate}
\end{thm}

Here, $\e(R)$ denotes the (Hilbert--Samuel) multiplicity of $\m$, while $K^R$ stands for the Koszul complex of a minimal system of generators of $\m$.

The conclusion of each of the three assertions of the above theorem says that the set
$$
\Delta(R)=\{t\in\R\mid\text{$\delta_t(G,X)\ne0$ for some generator $G\in\ds(R)$ and some object $X\in\ds(R)$}\}.
$$
of real numbers is bounded.
In the present paper, we shall look for classes of commutative noetherian rings $R$ such that the set $\Delta(R)$ is bounded.
The main result of this paper is the following theorem.

\begin{thm}\label{24}
Let $R$ be a commutative noetherian ring.
The set $\Delta(R)$ is bounded in the following cases.
\begin{enumerate}[\rm(1)]
\item
$R$ is a local ring with an isolated singularity.
\item
$R$ is a semilocal ring of Krull dimension at most one.
\item
$R$ is a Gorenstein semilocal J-0 domain of Krull dimension two.
\item
$R$ is a Gorenstein J-2 ring of finite Krull dimension.
\item
$R$ is a Cohen--Macaulay local J-0 domain of Krull dimension two, with a canonical module, and locally Gorenstein on the punctured spectrum.
\item
$R$ is a Cohen--Macaulay local J-2 ring, with a canonical module, and locally Gorenstein on the punctured spectrum.
\end{enumerate}
\end{thm}

Here, the J-0 and J-2 properties, together with the J-1 property, are classical notions concerning regular loci in prime spectra; see Definition \ref{16}.
By definition, every excellent ring is J-2.
Thus the following corollary is immediately deduced from (4) and (6) of the above theorem.

\begin{cor}
The set $\Delta(R)$ is bounded provided that $R$ is either
\begin{itemize}
\item
an excellent Gorenstein ring of finite Krull dimension, or
\item
an excellent Cohen--Macaulay local ring which admits a canonical module and which is locally Gorenstein on the punctured spectrum.
\end{itemize}
\end{cor}

In the process of proving the above theorem, we mainly consider how to generate syzygies from a single module.
We obtain the following theorem in this direction.

\begin{thm}\label{28}
Let $R$ be a commutative noetherian ring of Krull dimension $d$.
Let $R$ be one of the following:
$$
\begin{array}{lll}
&\text{{\rm(1)} a local ring with an isolated singularity;}\qquad
&\text{{\rm(2)} a semilocal ring with $d\le1$;}\qquad\\
&\text{{\rm(3)} a semilocal J-0 domain with $d=2$;}
&\text{{\rm(4)} a J-2 ring with $d<\infty$.}
\end{array}
$$
Then there exists a finitely generated $R$-module out of which the $d$th syzygy of each finitely generated $R$-module can be built only by taking direct summands and extensions finitely many times.
\end{thm}

This theorem proves under a weaker assumption a weaker version of the generation studied in \cite{radius,stgen,ua}, where the {\em strong generation} of the module category is discussed, which also evaluates the number of extensions necessary to build the high syzygies.
In relation to the assertion of the theorem in the cases of (1) and (4), the {\em classical generation} of the module category, which is weaker than our generation in that it allows taking not only direct summands and extensions, but also kernels of epimorphisms and cokernels of monomorphisms, is shown in \cite{kos} and \cite{jc}, respectively, the latter of which has been extended to schemes in \cite{ELS}.

The organization of this paper is as follows.
In Section 2, for a general triangulated category, we state the precise definition of a DHKK complexity and its basic properties.
We find out certain conditions for its non-vanishing range to have upper/lower bounds.
In Section 3, we show that generation of the high syzygies from a single module implies that the range is bounded below, and deal with the case of a semilocal ring of Krull dimension at most one.
In Section 4, we consider generating the maximal Cohen--Macaulay modules and the $d$th syzygies of finitely generated modules, where $d$ is the Krull dimension of the ring.
We then prove that the range is bounded under some mild assumptions.
Proofs of Theorems \ref{24} and \ref{28} are given at the end of this section.

%%%%%%%%%%%%%%%%%%%%%%%%%%%%%%%%%%%%%%%%%%
\section{DHKK complexities for a general triangulated category}

In this section, we give the definition of a DHKK complexity and its properties for a general triangulated category.
We also state and prove a key lemma at the end of this section.
First of all, we establish a setup.

\begin{setup}
Throughout this paper, all subcategories are assumed to be strictly full.
Throughout this section, let $\T$ be a triangulated category.
\end{setup}

We begin with recalling the definition of the binary operation $\star$ for subcategories of $\T$.

\begin{dfn}
For subcategories $\X,\Y$ of $\T$, we denote by $\X\star\Y$ the subcategory of $\T$ consisting of objects $T\in\T$ which fits into an exact triangle $X\to T\to Y\rightsquigarrow$ in $\T$ with $X\in\X$ and $Y\in\Y$.
\end{dfn}

As is shown in \cite[Lemma 2.3]{delta}, for subcategories $\X,\Y,\Z$ of $\T$ one has $(\X\star\Y)\star\Z=\X\star(\Y\star\Z)$.
Thus there is no danger of confusion even if for an integer $n\ge0$ and subcategories $\X,\X_1,\dots,\X_n$ of $\T$ we write
$$
\bigstar_{i=1}^n\X_i=\X_1\star\cdots\star\X_n,\qquad
\X^{\star n}=\underbrace{\X\star\cdots\star\X}_n.
$$
For an integer $n\ge0$ and objects $X,X_1,\dots,X_n$ of $\T$ we set $\bigstar_{i=1}^nX_i:=\bigstar_{i=1}^n\{X_i\}$ and $X^{\star n}:=\{X\}^{\star n}$.
The following lemma is also shown in \cite[Lemma 2.3]{delta}, which is used a couple of times later.

\begin{lem}\label{2}
Let $X_1,\dots,X_n$ be objects of $\T$.
Let $M$ be an object in $\bigstar_{i=1}^nX_i$.
Then $M[s]\in\bigstar_{i=1}^n(X_i[s])$ for all integers $s$, and $M\oplus(\bigoplus_{i=1}^nY_i)\in\bigstar_{i=1}^n(X_i\oplus Y_i)$ for all objects $Y_1,\dots,Y_n\in\T$.
\end{lem}

Now we recall the definition of a DHKK complexity, which has been introduced by Dimitrov, Haiden, Katzarkov and Kontsevich \cite{DHKK}.

\begin{dfn}[Dimitrov--Haiden--Katzarkov--Kontsevich]\label{1}
Let $X,Y\in\T$ and $t\in\R$.
We denote by $\delta_t(X,Y)$ the infimum of the sums $\sum_{i=1}^re^{n_it}$, where $r$ runs over the nonnegative integers and $n_i$ run over the integers such that there exists a sequence
$$
\xymatrix@R-2pc{
0\ar@{=}[r]& Y_0\ar[rr]&&Y_1\ar[r]\ar[dl]& \cdots\ar[r]& Y_{r-1}\ar[rr]&& Y_r\ar@{=}[r]\ar[dl]& Y\oplus Y'\\
&&X[n_1]\ar@{~>}[ul]&&\cdots&&X[n_r]\ar@{~>}[ul]
}
$$
of exact triangles $\{Y_{i-1}\to Y_i\to X[n_i]\rightsquigarrow\}_{i=1}^r$ in $\T$.
We call the function $\R\ni t\mapsto\delta_t(X,Y)\in\R_{\ge0}\cup\{\infty\}$ the {\em DHKK complexity} of $Y$ relative to $X$.
If $Y=0$, then one can take $r=0$, and hence $\delta_t(X,Y)=0$.
Also, by \cite[Proposition 2.6]{delta} one has the equality
$$
\textstyle
\delta_t(X,Y)=\inf\{\sum_{i=1}^r\e^{n_it}\mid\text{$r\in\ZZ_{\ge0}$ and $n_i\in\ZZ$ such that $Y\oplus Y'\in\bigstar_{i=1}^rX[n_i]$ for some $Y'\in\T$}\}.
$$
\end{dfn}

Next we recall several fundamental notions about triangulated categories.

\begin{dfn}
\begin{enumerate}[(1)]
\item
A {\em thick} subcategory of $\T$ is by definition a triangulated subcategory of $\T$ closed under direct summands.
\item
For an object $T$ of $\T$ we denote by $\thick_\T T$ the {\em thick closure} of $T$ in $\T$, that is, the smallest thick subcategory of $\T$ containing $T$.
\item
An object $G$ of $\T$ is called a {\em split generator} of $\T$ or a {\em classical generator} of $\T$ or a {\em thick generator} of $\T$, if $\thick_\T G=\T$.
In this paper we simply call such an object $G$ a {\em generator} of $\T$.
\item
We say that $\T$ is {\em periodic} if there exists a positive integer $n$ such that the $n$th shift functor $[n]$ is isomorphic to the identity functor of $\T$.
If $n$ is the least such integer, then $\T$ is called {\em $n$-periodic}.
\end{enumerate}
\end{dfn}

The following lemma, which is shown in \cite[Proposition 2.7 and Lemma 2.9]{delta}, provides a couple of properties of DHKK complexities.

\begin{lem}\label{3}
Let $X,Y,Z$ be objects of $\T$.
Let $t$ be a real number.
\begin{enumerate}[\rm(1)]
\item
One has $\delta_t(X,Y)<\infty$ if and only if $Y\in\thick_\T X$.
Hence $\delta_t(X,Y)$ is finite if $X$ is a generator of $\T$.
\item
Suppose that $\T$ is periodic and that $\delta_t(X,Y)$ is finite.
Then $\delta_t(X,Y)=0$ unless $t=0$.
\item
Suppose that both $\delta_t(X,Y)$ and $\delta_t(Y,Z)$ are finite.
One then has $\delta_t(X,Z)\le\delta_t(X,Y)\cdot\delta_t(Y,Z)$.
\end{enumerate}
\end{lem}

The lemma below is shown by applying the third assertion of the above lemma.

\begin{lem}\label{5}
Let $t$ be a real number such that $\delta_t(G,G)=0$ for some generator $G$ of $\T$.
Then $\delta_t(H,T)=0$ for every generator $H$ of $\T$ and every object $T$ of $\T$.
\end{lem}

\begin{proof}
The proof is similar to that of \cite[Lemma 3.7]{delta}.
By Lemma \ref{3}(1), the numbers $\delta_t(H,G)$ and $\delta_t(G,T)$ are finite.
By assumption, $\delta_t(G,G)$ vanishes, and in particular it is finite.
Applying Lemma \ref{3}(3), we have
$
0\le\delta_t(H,T)\le\delta_t(H,G)\cdot\delta(G,T)\le\delta_t(H,G)\cdot\delta_t(G,G)\cdot\delta_t(G,T)=0.
$
We obtain $\delta_t(H,T)=0$.
\end{proof}

Now we define the set $\Delta(\T)$ of real numbers which is the main target of this paper.

\begin{dfn}
We denote by $\Delta(\T)$ be the set of real numbers $t$ such that $\delta_t(G,X)\ne0$ for some generator $G$ of $\T$ and some object $X$ of $\T$.
\end{dfn}

We give here simple observations regarding the set $\Delta(\T)$.

\begin{prop}\label{25}
Suppose that $\T$ is not a zero category but admits a generator.
Then $\Delta(\T)\supseteq\{0\}$.
If $\T$ is moreover periodic, then $\Delta(\T)=\{0\}$.
\end{prop}

\begin{proof}
Let $G$ be a generator of $\T$.
Using Lemma \ref{3}(1), we observe that
$$
\delta_0(G,G)=\inf\{r\in\ZZ_{\ge0}\mid\text{$G\oplus G'\in\bigstar_{i=1}^rG[n_i]$ for some $G'\in\T$ and $n_i\in\ZZ$}\}\in\R_{\ge0}.
$$
As $\T$ is a nonzero category, $G$ is a nonzero object of $\T$.
Hence $r$ cannot be taken to be $0$, so that $\delta_0(G,G)\ge1$.
In fact, since $G\oplus0=G\in\bigstar_{i=1}^1G[0]$, we have $\delta_0(G,G)=1$.
Thus, $0$ belongs to $\Delta(\T)$.

If $\T$ is periodic, then it follows from Lemma \ref{3}(1)(2) that $\delta_t(C,T)=0$ for all generators $C$ of $\T$, all objects $T$ of $\T$, and all $0\ne t\in\R$, so that $\Delta(\T)\subseteq\{0\}$.
\end{proof}

Finally, we give a sufficient condition for $\Delta(\T)$ to be bounded above or below.
This lemma plays a key role in the proofs of our main results stated in the next two sections.

\begin{lem}\label{9}
Let $G$ be a generator of $\T$.
Let $n$ be an integer and $m$ a positive integer.
Suppose that $G[n]$ is isomorphic to a direct summand of some object in $G^{\star m}$.
Then the following hold.\\
\qquad{\rm(1)} If $n>0$, then $\sup\Delta(\T)\le\frac{1}{n}\log m$.
\qquad{\rm(2)} If $n<0$, then $\inf\Delta(\T)\ge\frac{1}{n}\log m$.
\end{lem}

\begin{proof}
By assumption, $G[n]\oplus H\in G^{\star m}$ for some $H\in\T$.
Putting $K=H[-n]$ and applying Lemma \ref{2}, we have that $G\oplus K\in(G[-n])^{\star m}$, that $G[-n]\oplus K[-n]\in(G[-2n])^{\star m}$, and that
$$
G\oplus K\oplus(K[-n])^{\oplus m}\in(G[-n]\oplus K[-n])^{\star m}\subseteq((G[-2n])^{\star m})^{\star m}=(G[-2n])^{\star m^2}.
$$
By a similar argument, we get
$$
G\oplus K\oplus(K[-n])^{\oplus m}\oplus(K[-2n])^{\oplus m^2}\in(G[-2n]\oplus K[-2n])^{\star m^2}\subseteq((G[-3n])^{\star m})^{\star m^2}=(G[-3n])^{\star m^3}.
$$
Iterating this procedure, for each integer $i>0$ we find an object $L_i\in\T$ such that $G\oplus L_i\in(G[-in])^{\star m^i}$.
Hence $0\le\delta_t(G,G)\le m^ie^{-int}=(me^{-nt})^i$.
Now let $t$ be such that $nt>\log m$.
Then $0<me^{-nt}<1$, so that $\lim_{i\to\infty}(me^{-nt})^i=0$, and we get $\delta_t(G,G)=0$.
It is seen from Lemma \ref{5} that $t\notin\Delta(\T)$.
Hence, $nt\le\log m$ for all $t\in\Delta(\T)$.
We conclude that $\sup\Delta(\T)\le\frac{1}{n}\log m$ if $n>0$ and $\inf\Delta(\T)\ge\frac{1}{n}\log m$ if $n<0$.
\end{proof}

%%%%%%%%%%%%%%%%%%%%%%%%%%%%%%%%%%%%%%
\section{Vanishing of DHKK complexities and generating syzygies}

From here to the end of this paper, we consider the boundedness of the non-vanishing range of DHKK complexities for the singularity category of a commutative noetherian ring $R$.
%In this section, we relate $\Delta(R)$ being bounded below to building the high syzygies out of a single module by taking direct summands and extensions.
We begin with a setup.

\begin{setup}
Throughout the rest of this paper, let $R$ be a commutative noetherian ring.
We abbreviate Krull dimension to dimension, and use the notation $\dim$.
We denote by $\mod R$ the category of finitely generated $R$-modules, by $\db(\mod R)$ the bounded derived category of the abelian category $\mod R$, and by $\ds(R)$ the {\em singularity category} of $R$, which is defined as the Verdier quotient of $\db(\mod R)$ by the thick closure $\thick_{\db(\mod R)}R$. 
The only triangulated category on which we work from here to the end of the paper is the singularity category $\ds(R)$.
So, for example, the operation $\star$ is always taken in $\ds(R)$.
We simply use the notation $\Delta(R)$ to denote $\Delta(\ds(R))$.
\end{setup}

We recall the definitions of an extension closure and a syzygy, both of which play central roles in the rest of this paper.

\begin{dfn}
\begin{enumerate}[(1)]
\item
A subcategory $\X$ of $\mod R$ is said to be {\em closed under extensions} provided that for each exact sequence $0\to L\to M\to N\to0$ in $\mod R$, if $L$ and $N$ belong to $\X$, then so does $M$.
For a finitely generated $R$-module $G$, we denote by $\ext G$ (or $\ext_RG$ to specify the base ring) the {\em extension closure} of $G$, which is defined as the smallest subcategory of $\mod R$ containing $G$ and closed under direct summands and extensions.
\item
Let $n\ge0$ be an integer.
Let $M$ be a finitely generated $R$-module.
We denote the {\em $n$th syzygy} of $M$ by $\syz^nM$ (or $\syz_R^nM$ to specify the base ring), which is defined by an exact sequence $0\to\syz^nM\to P_{n-1}\to\cdots\to P_1\to P_0\to M\to0$ in $\mod R$ with $P_i$ projective for all $0\le i\le n-1$.
Each finitely generated projective $R$-module $P$ is an $n$th syzygy for any $n\ge0$ since there is an exact sequence $0\to P\to P\to0\to\cdots\to0$.  
By Schanuel's lemma, $\syz^nM$ is uniquely determined by $M$ and $n$ up to projective summands.
There is an isomorphism $\syz^nM\cong M[-n]$ in $\ds(R)$ (see \cite[Lemma 2.4]{sing}).
We denote by $\syz^n(\mod R)$ the subcategory of $\mod R$ consiting of the $n$th syzygies of finitely generated $R$-modules.
\end{enumerate}
\end{dfn}

The first assertion of the lemma below relates a module in an extension closure to an object in a product with respect to the operation $\star$.
The second assertion says that a module is a generator of the singularity category if it generates all syzygies by direct summands and extensions.

\begin{lem}\label{7}
Let $G$ be a finitely generated $R$-module.
\begin{enumerate}[\rm(1)]
\item
For each $M\in\ext G$ there exist $N\in\mod R$ and $m\in\ZZ_{>0}$ such that $M\oplus N\in G^{\star m}$.
\item
If $\syz^n(\mod R)$ is contained in $\ext G$ for some nonnegative integer $n$, then $G$ is a generator of $\ds(R)$.
\end{enumerate}
\end{lem}

\begin{proof}
(1) One can deduce the assertion by making use of \cite[Proposition 2.4]{tes}.

(2) Fix an object $X\in\ds(R)$.
Then $X\cong M[r]\cong\syz^nM[n+r]$ in $\ds(R)$ for some $M\in\mod R$ and $r\in\ZZ$; see \cite[Lemma 2.4]{sing}.
By assumption, $\syz^nM$ belongs to $\ext G$, so that it is in $\thick_{\ds(R)}G$.
Hence $X$ is in $\thick_{\ds(R)}G$ as well, and we get $\thick_{\ds(R)}G=\ds(R)$.
Thus $G$ is a generator of $\ds(R)$.
\end{proof}

The following proposition gives a sufficient condition for $\Delta(R)$ to have a lower bound.

\begin{prop}\label{6}
Suppose that there exist a finitely generated $R$-module $G$ and a nonnegative integer $n$ such that $\syz^n(\mod R)$ is contained in $\ext G$.
Then the set $\Delta(R)$ is bounded below.
\end{prop}

\begin{proof}
Lemma \ref{7}(2) implies that $G$ is a generator of $\ds(R)$.
Since $\syz^{n+1}(\mod R)\subseteq\syz^n(\mod R)$, we may assume $n>0$.
By assumption, $\syz^nG$ belongs to $\ext G$.
Lemma \ref{7}(1) implies that $\syz^nG$ is isomorphic to a direct summand of an object in $G^{\star m}$ for some $m>0$.
Since $\syz^nG\cong G[-n]$ in $\ds(R)$ and $-n<0$, Lemma \ref{9}(2) implies that $\inf\Delta(R)\ge-\frac{1}{n}\log m$.
Thus $\Delta(R)$ is bounded below.
\end{proof}

Next we recall the classical notions of regular/singular loci and J-conditions.

\begin{dfn}\label{16}
We denote by $\reg R$ and $\sing R$ the {\em regular locus} and the {\em singular locus} of $R$, respectively.
That is, $\reg R$ is the set of prime ideals $\p$ of $R$ such that $R_\p$ is regular, and $\sing R$ is the complement of $\reg R$ in $\spec R$.
Following \cite[(32.B)]{M}, we say that $R$ is {\em J-0} if $\reg R$ contains a nonempty open subset of $\spec R$, that $R$ is {\em J-1} if $\reg R$ is an open subset of $\spec R$, and that $R$ is {\em J-2} if any finitely generated commutative $R$-algebra is J-1.
It follows from the definition that a residue ring of a J-2 ring is again J-2.
\end{dfn}

We need to state the lemma below, which should be a well-known fact.

\begin{lem}\label{8}
The following are equivalent.\\
\qquad
{\rm(1)} The set $\spec R$ is finite.\qquad
{\rm(2)} The ring $R$ is semilocal and has dimension at most one.\\
When one of these two conditions is satisfied, the ring $R$ is J-1.
\end{lem}

\begin{proof}
As there exist only finitely many minimal prime ideals, (2) implies (1).
If $R$ has dimension at least two, then there exist infinitely many prime ideals of $R$ with height one; see \cite[Theorem 144]{K}.
This shows that (1) implies (2).
Suppose that $\spec R$ is finite.
Then $\sing R$ is both specialization-closed and finite, so that it is closed.
Hence $\reg R$ is open.
Thus the second assertion follows.
\end{proof}

Using the above lemma, we get the result below.
The equality $\ext G=\mod R$ means that $G$ {\em generates} all finitely generated $R$-modules, or more precisely, that every finitely generated $R$-module can be built out of $G$ by taking direct summands and extensions finitely many times.

\begin{prop}\label{17}
Let $R$ be a semilocal ring of dimension at most one.
Then there exists a finitely generated $R$-module $G$ such that $\ext G=\mod R$.
\end{prop}

\begin{proof}
Put $G=\bigoplus_{\p\in\spec R}R/\p$.
Lemma \ref{8} says that $\spec R$ is a finite set.
Hence $G$ is a finitely generated $R$-module.
Recall that every finitely generated $R$-module has a filtration of submodules each of whose subquotients is isomorphic to $R/\p$ for some $\p\in\spec R$.
It is easy to observe that $\ext G=\mod R$.
\end{proof}

Now we can prove the following theorem, which is the main result of this section.

\begin{thm}\label{13}
Let $R$ be a semilocal ring of dimension at most one.
Then $\Delta(R)$ is bounded.
\end{thm}

\begin{proof}
Propositions \ref{6} and \ref{17} imply that $\Delta(R)$ is bounded below.
Let us show that it is bounded above.
Note that $\spec R=\Min R\cup\Max R$.
Put
$$
\textstyle
H=\bigoplus_{\p\in\Min R}R/\p,\qquad K=\bigoplus_{\m\in\Max R}R/\m,\qquad G=H\oplus K,\qquad C=\bigoplus_{\p\in\spec R}R/\p
$$
We have $G=C$ if $\dim R=1$, and $G=C^{\oplus2}$ if $\dim R=0$.
The $R$-modules $H,K,G,C$ are all finitely generated since $\spec R$ is finite by Lemma \ref{8}.
As is seen in the proof of Proposition \ref{17}, we have $\ext G=\ext C=\mod R$, so that $G$ is a generator of $\ds(R)$ by Lemma \ref{7}(2). 
Let $\p$ be a minimal prime ideal of $R$.
Then it is an associated prime ideal of $R$, and there is an injective homomorphism $R/\p\to R$.
There exists an exact sequence $0\to H\to F\to L\to0$ in $\mod R$ with $F$ free, and therefore $H\cong L[-1]$ in $\ds(R)$.
Let $\m$ be a maximal ideal of $R$.
Setting $g=\grade\m$ (now $g=0,1$ as $\dim R=1$), there is an $R$-regular sequence $\xx=x_1,\dots,x_g$ in $\m$.
Then $\depth R_\m=g$ by \cite[Proposition 1.2.10(a)]{BH}, and $\depth R_\m/\xx R_\m=0$.
Therefore $\m$ is an associated prime ideal of $R/(\xx)$, and we get an injective homomorphism $R/\m\to R/(\xx)$.
It is observed that there exists an exact sequence $0\to K\to P\to U\to0$ in $\mod R$ such that $P$ has finite projective dimension.
We have $K\cong U[-1]$ in $\ds(R)$.
It follows that $G=H\oplus K\cong Z[-1]$ in $\ds(R)$, where $Z:=L\oplus U$.
Since $Z\in\mod R=\ext G$, Lemma \ref{7}(1) implies that $Z$ is isomorphic to a direct summand of an object in $G^{\star m}$ for some $m>0$.
As $Z\cong G[1]$ in $\ds(R)$ and $1>0$, we see from Lemma \ref{9}(1) that $\Delta(R)$ is bounded above.
\end{proof}

Following the proofs of Theorem \ref{13} and Lemma \ref{9}, let us compute a concrete example.

\begin{ex}
Let $R=k[\![x,y,z]\!]/(x^2,y^2)$ where $k$ is a field.
Then $R$ is a local ring of dimension $1$.
We have $\Min R=\{\p\}$, $\Max R=\{\m\}$, and $\spec R=\{\p,\m\}$, where $\p=(x,y)$ and $\m=(x,y,z)$.
Set $G=R/\p\oplus R/\m$.

The element $z$ of $R$ is a non-zerodivisor, and $R/(z)$ has projective dimension $1$ over $R$.
There are exact sequences $0\to R/\p\xrightarrow{v}R\to R/(xy)\to0$ and $0\to R/\m\xrightarrow{w}R/(z)\to R/(xy,z)\to0$, where the maps $v,w$ are defined by $v(\overline1)=xy$ and $w(\overline1)=\overline{xy}$.
Taking the direct sum, we get an exact sequence $0\to G\to R\oplus R/(z)\to Z\to0$, where $Z=R/(xy)\oplus R/(xy,z)$.
Therefore $G\cong Z[-1]$ in $\ds(R)$.
There are descending chains $R/(xy)\supsetneq\p/(xy)\supsetneq(x)/(xy)\supsetneq0$ and $R/(xy,z)\supsetneq\m/(xy,z)\supsetneq(x,z)/(xy,z)\supsetneq0$.
It is seen that
$$
\begin{array}{l}
\p/(x)=(x,y)/(x)\cong R/(x):y=R/\p,\qquad
(x)/(xy)\cong R/(xy):x=R/\p,\\
\m/(x,z)=(x,y,z)/(x,z)\cong R/(x,z):y=R/\m,\qquad
(x,z)/(xy,z)\cong R/(xy,z):x=R/\m.
\end{array}
$$
Those two descending chains together with these four isomorphisms give rise to exact sequences
$$
\begin{array}{l}
0\to\p/(xy)\to R/(xy)\to R/\p\to0,\qquad
0\to R/\p\to\p/(xy)\to R/\p\to0,\\
0\to\m/(xy,z)\to R/(xy,z)\to R/\m\to0,\qquad
0\to R/\m\to\m/(xy,z)\to R/\m\to0.
\end{array}
$$
Taking the direct sums, we get exact sequences $0\to E\to Z\to G\to0$ and $0\to G\to E\to G\to0$, where $E=\p/(xy)\oplus\m/(xy,z)$.
It follows that $G[1]\cong Z$ belongs to $G^{\star3}$.
Using Lemma \ref{2}, we get $G\in(G[-1])^{\star3}$, which implies $G\in(G[-i])^{\star3^i}$ and $0\le\delta_t(G,G)\le3^ie^{-it}=(3e^{-t})^i$ for every integer $i>0$.
If $t>\log3$, then $0<3e^{-t}<1$ and taking $\lim_{i\to\infty}$ shows that $\delta_t(G,G)=0$.

As $G=R/\p\oplus R/\m$, we have $\syz G=\p\oplus\m$.
It is easy to observe that the following holds.
$$
\p\supsetneq(x)\cong R/0:x=R/(x)\supsetneq\p/(x)=(x,y)/(x)\cong R/(x):y=R/\p.
$$
We get exact sequences $0\to R/\p\to(x)\to R/\p\to0$ and $0\to(x)\to\p\to R/\p\to0$.
It is seen that $\p\in(R/\p)^{\star3}$.
Since $\m/\p=(x,y,z)/(x,y)\cong R/(x,y):z=R/\p$, there is an exact sequence $0\to\p\to\m\to R/\p\to0$.
Hence $\m\in(R/\p)^{\star3}\star(R/\p)=(R/\p)^{\star4}$.
Using the split exact sequence $0\to\p\to\p\oplus\m\to\m\to0$, we obtain $G[-1]\cong\syz G=\p\oplus\m\in(R/\p)^{\star3}\star(R/\p)^{\star4}=(R/\p)^{\star7}$.
Lemma \ref{2} implies $G[-1]\oplus(R/\m)^{\oplus7}\in(R/\p\oplus R/\m)^{\star7}=G^{\star7}$.
For each positive integer $i$ the object $G[-i]$ is isomorphic to a direct summand of some object in $G^{\star7^i}$, and $\delta_t(G,G)\le(7^ie^{it})=(7e^t)^i$.
We see that if $t<-\log7$, then $0<7e^t<1$ and $\delta_t(G,G)=0$.

Consequently, it holds that if $t>\log3$ or $t<-\log7$, then $\delta_t(D,X)=0$ for all generators $D$ and objects $X$ of $\ds(R)$ by Lemma \ref{5}.
This shows that one has the inclusion $\Delta(R)\subseteq[-\log7,\log3]$.
\end{ex}

%%%%%%%%%%%%%%%%%%%%%%%%%%%%
\section{Vanishing of DHKK complexities and generating maximal Cohen--Macaulay modules}

Continuing from the last section, we consider in this section the boundedness of the set $\Delta(R)$, paying attention to maximal Cohen--Macaulay modules over $R$.
First of all, let us recall the definition of a maximal Cohen--Macaulay module and related notions.

\begin{dfn}
\begin{enumerate}[(1)]
\item
A finitely generated $R$-module $M$ is called {\em maximal Cohen--Macaulay} if $\depth_{R_\p}M_\p\ge\dim R_\p$ for all prime ideals $\p$ of $R$.
Since $\depth0=\infty$, the zero module is maximal Cohen--Macaulay.
The ring $R$ is Cohen--Macaulay if and only if it is a maximal Cohen--Macaulay module over itself.
\item
We denote by $\cm(R)$ the subcategory of $\mod R$ consisting of maximal Cohen--Macaulay $R$-modules.
Note that $\cm(R)$ contains $0$ and is closed under direct summands and extensions.
\item
Suppose that $R$ is a Cohen--Macaulay ring.
A finitely generated $R$-module $\omega$ is a {\em canonical module} of $R$ if for each prime ideal $\p$ of $R$ the localization $\omega_\p$ is a canonical module of the local ring $R_\p$; see \cite[Definition 3.3.16]{BH}.
When the ring $R$ is Gorenstein, the $R$-module $R$ is a canonical module of $R$.
\end{enumerate}
\end{dfn}

From now on we state and prove two lemmas.
The first one says that taking the syzygy and the canonical dual are compatible with taking the extension closure.

\begin{lem}\label{22}
Let $G$ and $M$ be finitely generated $R$-modules such that $M\in\ext G$.
\begin{enumerate}[\rm(1)]
\item
For every nonnegative integer $n$ it holds that $\syz^nM\in\ext(R\oplus\syz^nG)$.
\item
Suppose that $R$ is a Cohen--Macaulay ring with a canonical module $\omega$ and $G$ is maximal Cohen--Macaulay.
It then holds that $M^\dag\in\ext(G^\dag)$, where $(-)^\dag=\Hom_R(-,\omega)$.
\end{enumerate}
\end{lem}

\begin{proof}
(1) Let $\X$ be the subcategory of $\mod R$ consisting of modules $X$ such that $\syz^nX\in\ext(R\oplus\syz^nG)$.
If $Y$ is a direct summand of $X\in\X$, then $\syz^nY$ is a dierct summand of $\syz^nX\in\ext(R\oplus\syz^nG)$, and hence $\syz^nY\in\ext(R\oplus\syz^nG)$, which means $Y\in\X$.
Let $0\to A\to B\to C\to0$ be an exact sequence in $\mod R$ with $A,C\in\X$.
Then there exists an exact sequence $0\to\syz^nA\to\syz^nB\to\syz^nC\to0$ in $\mod R$ up to projective summands, and $\syz^nA,\syz^nC$ are in $\ext(R\oplus\syz^nG)$.
Hence $\syz^nB$ belongs to $\ext(R\oplus\syz^nG)$, which means that $B$ belongs to $\X$.
Thus $\X$ is closed under direct summands and extensions.
Since $\X$ contains $G$, it also contains $\ext G$.
This particularly says that the module $M$ belongs to $\X$, which means that the assertion holds.

(2) Let $\X$ be the subcategory of $\mod R$ consisting of maximal Cohen--Macaulay $R$-modules $X$ such that $X^\dag\in\ext(G^\dag)$.
If $Y$ is a direct summand of $X\in\X$, then $Y^\dag$ is a dierct summand of $X^\dag\in\ext(G^\dag)$, and hence $Y^\dag\in\ext(G^\dag)$, which means $Y\in\X$.
Let $0\to A\to B\to C\to0$ be an exact sequence in $\mod R$ with $A,C\in\X$.
Then $A,C$ are maximal Cohen--Macaulay $R$-modules, and so is $B$.
Applying $(-)^\dag$ gives rise to an exact sequence $0\to C^\dag\to B^\dag\to A^\dag\to0$ of maximal Cohen--Macaulay $R$-modules, since $\Ext_R^1(C,\omega)=0$.
Since $A^\dag,C^\dag$ are in $\ext(G^\dag)$, so is $B^\dag$, which means that $B$ belongs to $\X$.
Thus $\X$ is closed under direct summands and extensions.
Since $\X$ contains $G$, it also contains $\ext G$.
This particularly says that the module $M$ belongs to $\X$, which means that the assertion holds.
\end{proof}

In the lemma below, we study the condition that an extension closure contains all the maximal Cohen--Macaulay modules.
The second assertion is similar to Lemma \ref{7}(2).
The {\em (large) restricted flat dimension} $\rfd_RM$ of a finitely generated $R$-module $M$ is defined by $\rfd_RM=\sup_{\p\in\spec R}(\depth R_\p-\depth M_\p)$.

\begin{lem}\label{20}
Assume that the ring $R$ is Cohen--Macaulay.
Let $G$ be a finitely generated $R$-module.
\begin{enumerate}[\rm(1)]
\item
There exists an integer $n\ge0$ such that $\syz^nG$ is a maximal Cohen--Macaulay $R$-module.
\item
If $\cm(R)$ is contained in $\ext G$, then $G$ is a generator of $\ds(R)$.
\item
Suppose that $R$ has a canonical module $\omega$.
Put $(-)^\dag=\Hom_R(-,\omega)$.
The following two statements hold.
\begin{enumerate}[\rm(a)]
\item
If $G$ is maximal Cohen--Macaulay, then for any $n\ge0$ there is an exact sequence $0\to G\to\omega^{\oplus e_0}\to\omega^{\oplus e_1}\to\cdots\to\omega^{\oplus e_{n-1}}\to H\to0$ of maximal Cohen--Macaulay $R$-modules, where $H=(\syz^n(G^\dag))^\dag$.
\item
If $\cm(R)\subseteq\ext G$, then there exists a maximal Cohen--Macaulay $R$-module $C$ with $\cm(R)=\ext C$.
\end{enumerate}
\end{enumerate}
\end{lem}

\begin{proof}
(1) It follows from \cite[Theorem 1.1]{AIL} that $r:=\rfd_RG$ is a nonnegative integer.
Take any integer $n\ge r$.
Then it is easy to observe that $\syz^nG$ is a maximal Cohen--Macaulay $R$-module.

(2) Let $X$ be any object of $\ds(R)$.
Then there exist a finitely generated $R$-module $N$ and an integer $n$ such that $X\cong N[n]$ in $\ds(R)$ by \cite[Lemma 2.4]{sing}.
According to (1), we can choose an integer $r\ge0$ such that $L:=\syz^rN$ is maximal Cohen--Macaulay.
Then $X\cong L[r+n]$ in $\ds(R)$.
As $L$ is in $\cm(R)$, which is contained in $\ext G$, we see that $L$ belongs to $\thick_{\ds(R)}G$ and so does $X$.
Thus $\thick_{\ds(R)}G=\ds(R)$.

(3a) Take an exact sequence $0\to\syz^n(G^\dag)\to R^{\oplus e_{n-1}}\to\cdots\to R^{\oplus e_1}\to R^{\oplus e_0}\to G^\dag\to0$ in $\mod R$.
Note that all the terms in this exact sequence belong to $\cm(R)$.
In general, the natural map $X\to X^{\dag\dag}$ is an isomorphism and $\Ext_R^i(X,\omega)=0$ for all $X\in\cm(R)$ and $i>0$.
Applying $(-)^\dag$ yields an exact sequence $0\to G\to \omega^{\oplus e_0}\to\cdots\to\omega^{\oplus e_{n-1}}\to H\to0$, where $H:=(\syz^n(G^\dag))^\dag\in\cm(R)$.

(3b) By (1) the $R$-module $\syz^nG$ is maximal Cohen--Macaulay for some $n\ge0$.
Fix $X\in\cm(R)$.
By (3a) (and its proof) there exist exact sequences of maximal Cohen--Macaulay $R$-modules:
$$
\begin{array}{l}
0\to\syz^nG\to \omega^{\oplus e_0}\to\cdots\to\omega^{\oplus e_{n-1}}\to(\syz^n((\syz^nG)^\dag))^\dag\to0,\\
0\to X^\dag\to\omega^{\oplus c_0}\to\cdots\to\omega^{\oplus c_{n-1}}\to(\syz^nX)^\dag\to0.
\end{array}
$$
It follows from \cite[Lemma 5.8]{radius} that $X^\dag\in\ext(R\oplus\syz^n((\syz^nX)^\dag)\oplus(\bigoplus_{i=0}^{n-1}\syz^i\omega))$.
Lemma \ref{22}(2) implies that
\begin{equation}\label{26}
\textstyle
X\cong X^{\dag\dag}\in\ext(\omega\oplus(\syz^n((\syz^nX)^\dag))^\dag\oplus(\bigoplus_{i=0}^{n-1}(\syz^i\omega)^\dag)).
\end{equation}
Note that $X$ is in $\cm(R)$, which is contained in $\ext G$, whence $X\in\ext G$.
Applying Lemma \ref{22} repeatedly, we get $\syz^nX\in\ext(R\oplus\syz^nG)$, $(\syz^nX)^\dag\in\ext(\omega\oplus(\syz^nG)^\dag)$, $\syz^n((\syz^nX)^\dag)\in\ext(R\oplus\syz^n\omega\oplus\syz^n((\syz^nG)^\dag))$, and
\begin{equation}\label{27} (\syz^n((\syz^nX)^\dag))^\dag\in\ext(\omega\oplus(\syz^n\omega)^\dag\oplus(\syz^n((\syz^nG)^\dag))^\dag).
\end{equation}
Combining \eqref{26} and \eqref{27}, we see that $X$ is in $\ext C$, where $C:=\omega\oplus(\bigoplus_{i=0}^n(\syz^i\omega)^\dag)\oplus(\syz^n((\syz^nG)^\dag))^\dag\in\cm(R)$.
Note that $C$ is independent of the choice of $X$.
We obtain $\cm(R)=\ext C$.
\end{proof}

Now we can prove the following proposition, whose second assertion gives a sufficient condition for $\Delta(R)$ to be bounded.
Recall that the {\em punctured spectrum} of a local ring $(R,\m)$ is by definition the set $\spec R\setminus\{\m\}$.

\begin{prop}\label{11}
Suppose either that $R$ is a Gorenstein ring, or that $R$ is a Cohen--Macaulay local ring which has a canonical module and is locally Gorenstein on the punctured spectrum.
\begin{enumerate}[\rm(1)]
\item
For each $M\in\cm(R)$, there exists $N\in\cm(R)$ such that $M\in\ext(R\oplus\syz N)$.
\item
If $\syz^n(\mod R)\subseteq\ext G$ for some $n\ge0$ and some finitely generated $R$-module $G$, then $\Delta(R)$ is bounded.
\end{enumerate}
\end{prop}

\begin{proof}
(1) Let $R$ be a Cohen--Macaulay ring.
Let $\omega$ be a canonical module of $R$; we take $\omega=R$ when $R$ is Gorenstein.
Lemma \ref{20}(3a) gives an exact sequence $0\to M\to\omega^{\oplus a}\xrightarrow{v}C\to0$ of maximal Cohen--Macaulay $R$-modules.
Take an exact sequence $0\to\syz C\to R^{\oplus b}\xrightarrow{w}C\to0$.
The pullback diagram of $v$ and $w$ induces an exact sequence $0\to\syz C\to M\oplus R^{\oplus b}\to\omega^{\oplus a}\to0$.
If $R$ is Gorenstein, then $\omega=R$ and $M\in\ext(R\oplus\syz C)$, which shows the assertion.
So, assume that $(R,\m,k)$ is a $d$-dimensional local ring which is locally Gorenstein on the punctured spectrum.
Then $\omega$ is locally free on the punctured spectrum of $R$, and hence it belongs to $\ext(\syz^dk)$ by \cite[Corollary 2.6]{stcm}.
It follows from \cite[Theorem 4.1(2)]{restf} that $\syz^dk$ is the $(d+1)$st syzygy of some finitely generated $R$-module $L$.
Setting $N=C\oplus\syz^dL\in\cm(R)$, we have $M\in\ext(R\oplus\syz N)$ and are done.

(2) Proposition \ref{6} implies that $\Delta(R)$ is bounded below.
Let $\omega$ be a canonical module of $R$.
Fix $M\in\cm(R)$.
By Lemma \ref{20}(3a) and \cite[Lemma 5.8]{radius} there exists a (maximal Cohen--Macaulay) $R$-module $N$ such that $M$ is in $\ext(R\oplus\syz^nN\oplus(\bigoplus_{i=0}^{n-1}\syz^i\omega))$.
Setting $H=R\oplus G\oplus(\bigoplus_{i=0}^{n-1}\syz^i\omega)$, we see that $\cm(R)$ is contained in $\ext H$.
In view of Lemma \ref{20}(3b), we may assume that $H$ is maximal Cohen--Macaulay and $\cm(R)=\ext H$.
Applying assertion (1) to $H$, we find a maximal Cohen--Macaulay $R$-module $K$ such that $H\in\ext(R\oplus\syz K)$.
Then $\cm(R)=\ext H\subseteq\ext(R\oplus\syz K)\subseteq\cm(R)$, and we get $\cm(R)=\ext(R\oplus\syz K)$.
Note that $R\oplus\syz K\cong K[-1]$ in $\ds(R)$.
As $K\in\cm(R)=\ext(R\oplus\syz K)$, by Lemma \ref{7}(1) we find a finitely generated $R$-module $L$ and an integer $m>0$ such that $K\oplus L\in(R\oplus\syz K)^{\star m}=(K[-1])^{\star m}$.
Using Lemma \ref{2}, we get $K[1]\oplus L[1]\in K^{\star m}$.
Lemma \ref{20}(2) implies that $R\oplus\syz K$ is a generator of $\ds(R)$, and so is $K$.
It follows from Lemma \ref{9}(1) that $\Delta(R)$ is bounded above.
Now the proof of the assertion is completed.
\end{proof}

%\begin{lem}\label{10}
%Let $R$ be a Gorenstein ring of dimension $d<\infty$. Suppose that there exist an integer $n\ge0$ and a finitely generated $R$-module $G$ such that $\syz^n(\mod R)\subseteq\ext G$. Then there exists a maximal Cohen--Macaulay $R$-module $C$ such that $\cm(R)=\ext C$.
%\end{lem}

%\begin{proof}
%It holds that $\cm(R)=\syz^{d+n}(\mod R)\subseteq\ext(R\oplus\syz^dG)\subseteq\cm(R)$. Indeed, the equality comes from the assumption that $R$ is Gorenstein, the first inclusion follows by Lemma \ref{22}(1), and the last inclusion is obvious. Setting $C=R\oplus\syz^dG$, we get $\cm(R)=\ext C$.
%\end{proof}

To apply the second assertion of the above proposition to deduce the boundedness of $\Delta(R)$, we need to consider when $\syz^n(\mod R)\subseteq\ext G$ for some $n\ge0$.
For this, we establish a lemma.

\begin{lem}\label{18}
Put $d=\dim R$ and assume $0<d<\infty$.
Let $x$ be a non-zerodivisor of $R$.
\begin{enumerate}[\rm(1)]
\item
One has the inequality $\dim R/(x)\le d-1$.
\item
Assume the localization $R_x$ is a regular ring.
Suppose there exist an finitely generated $R/(x)$-module $G$ and a positive integer $n$ such that $\syz_{R/(x)}^{d-1}(\mod R/(x))\subseteq\ext_{R/(x)}G$.
Then $\syz_R^d(\mod R)\subseteq\ext_R(R\oplus\syz_RG)$.
\end{enumerate}
\end{lem}

\begin{proof}
(1) Assume contrarily that $\dim R/(x)\ge d$.
Then there exists an ascending chain $(x)\subseteq\p_0\subsetneq\p_1\subsetneq\cdots\subsetneq\p_d$ of ideals of $R$ such that each $\p_i$ is prime.
Since $x$ is a non-zerodivisor, $\p_0$ has height at least one, so that it strictly contains some prime ideal $\q$ of $R$.
We get a chain $\q\subsetneq\p_0\subsetneq\cdots\subsetneq\p_d$ in $\spec R$, which implies that $\dim R\ge d+1$.
This contradicts the fact that $\dim R=d$.
Thus we obtain $\dim R/(x)\le d-1$.

(2) Fix a finitely generated $R$-module $M$.
Set $N=\syz^dM$.
The global dimension of the ring $R_x$ is equal to $\dim R_x$, which is at most $\dim R=d$.
Hence $\Ext_R^1(N,\syz N)_x\cong\Ext_R^{d+1}(M,\syz N)_x\cong\Ext_{R_x}^{d+1}(M_x,(\syz N)_x)=0$, and therefore $x^r\Ext_R^1(N,\syz N)=0$ for some positive integer $r$.
The module $N=\syz^dM$ is torsion-free as $d>0$, so that $x^r$ is a non-zerodivisor on $N$.
By \cite[Remark 2.12]{ua}, the module $N$ is isomorphic to a direct summand of $\syz_R(N/x^rN)$.
There is a series $\{0\to N/xN\xrightarrow{x^i}N/x^{i+1}N\to N/x^iN\to0\}_{i=1}^{r-1}$ of exact sequences in $\mod R$, which shows that $N/x^rN$ belongs to $\ext_R(N/xN)$.
Lemma \ref{22}(1) implies that $\syz_R(N/x^rN)$ belongs to $\ext_R(R\oplus\syz_R(N/xN))$, and so does $N$.
There is a commutative diagram
$$
\xymatrix@R-1pc@C3pc{
0\ar[r]& \syz^dM\ar[r]\ar[d]^x& P_{d-1}\ar[r]\ar[d]^x& \cdots\ar[r]& P_1\ar[r]\ar[d]^x& \syz M\ar[r]\ar[d]^x& 0\\
0\ar[r]& \syz^dM\ar[r]& P_{d-1}\ar[r]& \cdots\ar[r]& P_1\ar[r]& \syz M\ar[r]& 0
}
$$
of exact sequences in $\mod R$ with each $P_i$ projective.
The vertical arrows are the multiplication maps by $x$, all of which are injective.
The snake lemma induces an exact sequence
$$
0\to\syz^dM/x\syz^dM\to P_{d-1}/xP_{d-1}\to\cdots\to P_1/xP_1\to\syz M/x\syz M\to0
$$
in $\mod R/(x)$, which gives rise to an isomorphism $N/xN=\syz_R^dM/x\syz_R^dM\cong\syz_{R/(x)}^{d-1}(\syz_RM/x\syz_RM)$ up to $R/(x)$-projective summands.
By assumption, $N/xN$ belongs to $\ext_{R/(x)}G$.
We observe that $N/xN$ belongs to $\ext_RG$ as well.
Applying Lemma \ref{22}(1) again, we get $\syz_R(N/xN)\in\ext_R(R\oplus\syz_RG)$.
Consequently, we obtain $N\in\ext_R(R\oplus\syz_R(N/xN))\subseteq\ext_R(R\oplus\syz_RG)$.
Thus $\syz_R^d(\mod R)$ is contained in $\ext_R(R\oplus\syz_RG)$
\end{proof}

Let $M$ be a finitely generated $R$-module.
The {\em nonfree locus} $\nf(M)$ of $M$ is by definition the set of prime ideals $\p$ of $R$ such that $M_\p$ is nonfree as an $R_\p$-module.
For a sequence $\xx=x_1,\dots,x_n$ of elements of $R$ and an integer $1\le i\le n$ we denote by $\h_i(\xx,M)$ the $i$th {\em Koszul homology} of $\xx$ with coefficients in $M$.
A local ring $R$ is said to be with an {\em isolated singularity} if $R$ is locally regular on the punctured spectrum, or equivalently, if $\sing R\subseteq\{\m\}$ where $\m$ is the maximal ideal of $R$.
Using the above lemma, we can show the proposition below.
This result provides a certain class of rings that admit a module which generates all high syzygies. 

\begin{prop}\label{12}
Put $d=\dim R$.
There exists a finitely generated $R$-module $G$ such that $\syz^d(\mod R)$ is contained in $\ext G$, if either {\rm(1)} $R$ is a local ring with an isolated singularity, or
{\rm(2)} $R$ is a semilocal J-0 domain with $d=2$, or
{\rm(3)} $R$ is a J-2 ring with $d<\infty$.
\end{prop}

\begin{proof}
(1) Let $\m$ be the maximal ideal of the local ring $R$ and $k$ the residue field of $R$.
Fix $M\in\syz^d(\mod R)$.
As $R$ has an isolated singularity, we have that $\nf(M)\subseteq\{\m\}=\V(\m)$.
According to \cite[Lemma 3.4]{kos}, we can choose a system of parameters $\xx=x_1,\dots,x_d$ of $R$ such that $\nf(M)\subseteq\V(\xx)$ and $\xx\Ext_R^i(M,N)=0$ for all integers $i>0$ and all $R$-modules $N$.
It follows from \cite[Corollary 3.2(2)]{kos} that there exists a series $\{0\to\h_i(\xx,M)\to E_i\to\syz E_{i-1}\to0\}_{i=1}^d$ of exact sequences in $\mod R$ such that $E_0=\h_0(\xx,M)$ and $M$ is a direct summand of $E_d$.
For each integer $0\le i\le d$, the $R$-module $\h_i(\xx,M)$ has finite length since $\xx\h_i(\xx,M)=0$, and hence it belongs to the extension closure $\ext_Rk$.
An inductive argument shows that $M$ belongs to $\ext G$, where $G:=R\oplus(\bigoplus_{i=0}^{d-1}\syz^ik)\in\mod R$.
Consequently, $\syz^d(\mod R)$ is contained in $\ext G$.

(2) Since $R$ is J-0, there exists an element $x\in R$ such that $\emptyset\ne\D(x)\subseteq\reg R$.
Then $x$ is nonzero (hence it is a non-zerodivisor as $R$ is a domain), and the ring $R_x$ is regular.
Lemma \ref{18}(1) says $\dim R/(x)\le1$.
As $R/(x)$ is semilocal, Proposition \ref{17} gives rise to a finitely generated $R/(x)$-module $C$ such that $\ext_{R/(x)}C=\mod R/(x)$.
Therefore, we have that $\syz_{R/(x)}^{2-1}(\mod R/(x))\subseteq\mod R/(x)=\ext_{R/(x)}C$.
It follows from Lemma \ref{18}(2) that $\syz_R^2(\mod R)$ is contained in $\ext_R(R\oplus\syz_RC)$.
Setting $G=R\oplus\syz_RC$, we are done.

(3) We prove the proposition by induction on $d$.
We begin with handling the case $d=0$.
In this case, the ring $R$ is artinian and hence it is semilocal.
By Proposition \ref{17} we find a finitely generated $R$-module $G$ such that $\ext G=\mod R=\syz^d(\mod R)$ (as $d=0$).
Thus the assertion follows.
So, we assume that $d>0$.

(a) First we consider the case where $R$ is a domain.
Since $R$ is J-2, it is J-1 and hence $\reg R=\D(I)$ for some ideal $I$ of $R$.
As $R$ is a domain, the zero ideal of $R$ belongs to $\reg R$, so that we find a nonzero element $x\in I$.
Then $\D(x)\subseteq\D(I)=\reg R$, which implies that the localization $R_x$ is a regular ring.
Since $R$ is J-2, so is the residue ring $R/(x)$.
Lemma \ref{18}(1) says that $n:=\dim R/(x)\le d-1$.
The induction hypothesis yields a finitely generated $R/(x)$-module $C$ such that $\syz_{R/(x)}^{d-1}(\mod R/(x))\subseteq\syz_{R/(x)}^n(\mod R/(x))\subseteq\ext_{R/(x)}C$.
Putting $G=R\oplus\syz_RC$, we see from Lemma \ref{18}(2) that $\syz_R^d(\mod R)$ is contained in $\ext_RG$.

(b) Next we consider the general case.
There exists a chain $0=I_0\subsetneq\cdots\subsetneq I_m=R$ of ideals of $R$ such that for each $1\le i\le m$ one has $I_i/I_{i-1}\cong R/\p_i$ where $\p_i\in\spec R$.
As $R$ is J-2, so is the residue ring $R/\p_i$.
We have $d_i:=\dim R/\p_i\le d$.
Since $R/\p_i$ is a domain, by (a) we find a finitely generated $R/\p_i$-module $G_i$ such that $\syz_{R/\p_i}^d(\mod R/\p_i)\subseteq\syz_{R/\p_i}^{d_i}(\mod R/\p_i)\subseteq\ext_{R/\p_i}G_i$.
Put $G=R\oplus(\bigoplus_{i=1}^m(R/\p_i\oplus G_i\oplus(\bigoplus_{j=0}^{d-1}\syz_R^j\p_i)))\in\mod R$.

Fix an integer $1\le i\le m$ and a finitely generated $R/\p_i$-module $M$.
We have $\syz_{R/\p_i}^dM\in\syz_{R/\p_i}^d(\mod R/\p_i)\subseteq\ext_{R/\p_i}G_i$.
By \cite[Proposition 5.3]{radius} we see that $\syz_R^dM\in\ext_R(R\oplus R/\p_i\oplus\syz_{R/\p_i}^dM\oplus(\bigoplus_{j=0}^{d-1}\syz_R^j\p_i))\subseteq\ext_RG$.

Now, fix a finitely generated $R$-module $N$.
For each integer $1\le i\le m$ there exists an exact sequence $0\to I_{i-1}N\to I_iN\to I_iN/I_{i-1}N\to0$ in $\mod R$, which induces an exact sequence $0\to\syz_R^d(I_{i-1}N)\to\syz_R^d(I_iN)\to\syz_R^d(I_iN/I_{i-1}N)\to0$ in $\mod R$ up to projective summands.
Since $I_iN/I_{i-1}N$ is an $R/\p_i$-module, we see that $\syz_R^d(I_iN/I_{i-1}N)\in\ext_RG$.
An inductive argument shows that $\syz_R^dN$ belongs to $\ext_RG$.
Thus we conclude that $\syz_R^d(\mod R)$ is contained in $\ext_RG$.
\end{proof}

As a direct consequence of the combination of Propositions \ref{11}(2) and \ref{12}, we obtain the following theorem.

\begin{thm}\label{14}
Let $R$ be either a local ring with an isolated singularity, or a semilocal J-0 domain of dimension two, or a J-2 ring of finite dimension.
Suppose that $R$ is either Gorenstein, or a Cohen--Macaulay local ring with a canonical module and locally Gorenstein on the punctured spectrum.
Then the set $\Delta(R)$ is bounded.
\end{thm}

Finally, we give proofs of the theorems stated in the Introduction (Section 1).

\begin{proof}[Proof of Theorem \ref{24}]
The assertion in cases (2)--(6) follows from Theorems \ref{13} and \ref{14}.
The assertion in case (1) is a direct consequence of Theorem \ref{23}(1).
\end{proof}

\begin{proof}[Proof of Theorem \ref{28}]
The assertion is shown by Propositions \ref{17}(2) and \ref{12}.
\end{proof}

\begin{rem}
A direct consequence of Propositions \ref{11}(2) and \ref{12}(1) is that $\Delta(R)$ is bounded if $R$ is a {\em Cohen--Macaulay} local ring with an isolated singularity and {\em admitting a canonical module}.
Theorem \ref{23}(1) says that $\Delta(R)$ is bounded whenever $R$ is a local ring with an isolated singularity.
\end{rem}

%%%%%%%%%%%%%%%%%%%%%%%%%%%%%%%%%%%%%%

\end{document}